
\magnification=1050
\def\arXiv{oui}
  
 \def\oui{oui} 
  
\ifx\arXiv\oui
\else
 \pdfpagewidth=210truemm
 \pdfpageheight=297truemm 
\fi
  
  %
  %

  %
  %

  \catcode`@=12 

 \def\defrefnote#1{\definexref{#1}{{\the\footnotenumber}}{refnotes}}

  %
  %


\ifx\couleurs\oui
\input graphicx
 \pdfpagewidth=210truemm
 \pdfpageheight=297truemm 
 \voffset=-5mm
\fi

\input eplain.tex
\expandafter\def\expandafter\newdimen\expandafter{\newdimen}

\ifx\couleurs\oui
\beginpackages
\usepackage{color}
\endpackages 
 \pdfpagewidth=210truemm
 \pdfpageheight=297truemm 
\long\def\rge#1{{\color{red}#1}}

\definecolor{bleu-iecn}{cmyk}{.98,.13,.1,.55}

\else
\long\def\rge#1{#1}

\fi

\makeatletter
\def\numberedfootnote{%
ÊÊ\global\advance\footnotenumber by 1
ÊÊ\@eplainfootnote{{\number\footnotenumber}}%
}%
\def\makecolumns#1/#2 {\par \begingroup
ÊÊ \@columndepth = #1
ÊÊ \advance\@columndepth by -1
ÊÊ \divide \@columndepth by #2
ÊÊ \advance\@columndepth by 1
ÊÊ \@linestogoincolumn = \@columndepth
ÊÊ \@linestogo = #1
ÊÊ \currentcolumn = 1
ÊÊ \def\@endcolumnactions{%
ÊÊÊÊÊÊ\ifnum \@linestogo<2
ÊÊÊÊÊÊÊÊ \the\crtok \egroup \endgroup \par 
ÊÊÊÊÊÊ\else
ÊÊÊÊÊÊÊÊ \global\advance\@linestogo by -1
ÊÊÊÊÊÊÊÊ \ifnum\@linestogoincolumn<2
ÊÊÊÊÊÊÊÊÊÊÊÊ\global\advance\currentcolumn by 1
ÊÊÊÊÊÊÊÊÊÊÊÊ\global\@linestogoincolumn = \@columndepth
ÊÊÊÊÊÊÊÊÊÊÊÊ\the\crtok
ÊÊÊÊÊÊÊÊ \else
ÊÊÊÊÊÊÊÊÊÊÊÊ&\global\advance\@linestogoincolumn by -1
ÊÊÊÊÊÊÊÊ \fi
ÊÊÊÊÊÊ\fi
ÊÊ }%
ÊÊ \makeactive\^^M
ÊÊ \letreturn \@endcolumnactions
ÊÊ \@columnwidth = \hsize
ÊÊÊÊ \advance\@columnwidth by -\parindent
ÊÊÊÊ \divide\@columnwidth by #2
ÊÊ \penalty\abovecolumnspenalty
ÊÊ \noindent 
ÊÊ \valign\bgroup
ÊÊÊÊ &\hbox to \@columnwidth{\strut \hsize = \@columnwidth ##\hfil}\cr
}%
\makeatother

\lefteqnumbers
   \def\testd{oui}
   \def\choixlat{\ifx\numadroite\testd\righteqnumbers
            \else  \lefteqnumbers\fi}
    \choixlat

\catcode`@=\letter
\def\@eplainfootnote#1{\let\@sf\empty 
  \ifhmode\edef\@sf{\spacefactor\the\spacefactor}\/\fi
  \global\advance\hlfootlabelnumber by 1
  \hlstart@impl{foot}{\hlfootlabel}%
  \hldest@impl{footback}{\hlfootbacklabel}%
  \hbox{$^{(#1)}$}%
  \hlend@impl{foot}%
  \@sf\vfootnote{#1.}%
}%
\catcode`@=\other

  \interfootnoteskip=0pt
  \let\note=\numberedfootnote
  \everyfootnote={\eightpoint\leftskip=5truemm\rightskip5truemm}
  
  \hsize150truemm\vsize 240truemm\hoffset=5truemm

  \pretolerance=500\tolerance=1000\brokenpenalty=5000
  \parindent3mm
  
  \countdef\temps=170
  \temps=\time
  \countdef\nminutes=171{\nminutes = \time}
  \countdef\nheure=172
  \def\heure{\begingroup                     
     \temps = \time \divide\temps by 60
     \nheure = \temps                        
     \nminutes = \time
     \multiply\temps by 60
     \advance\nminutes by -\temps            
     \ifnum\nminutes<10 \toks1 = {0}%
     \else\toks1 = {}%
     \fi
     \number\nheure h\the\toks1 \number\nminutes  
  \endgroup}%

  \newcount\chstart
  \chstart=\pageno
 \headline={\ifnum\pageno=\chstart {\hfill} \else {\hss \tenrm --\ \folio\ --\hss}\fi}
  \footline={\hfill}
  \normalbaselines
  \frenchspacing
    \def\dater{\vglue-10mm\rightline{(\the\day/\the\month/\the\year)}}
  \def\dateheure{\vglue-10mm\rightline{(\the\day/\the\month/\the\year,\ \heure)}}

  \newif\ifpagetitre \pagetitretrue
\newtoks\hautpagetitre \hautpagetitre={\hfill}
\newtoks\baspagetitre \baspagetitre={\hfill}
\newtoks\auteurcourant \auteurcourant={\hfill}
\newtoks\titrecourant \titrecourant={\hfill}
\newtoks\hautpagegauche
\newtoks\hautpagedroite
\newtoks\hautpagemilieu
\hautpagemilieu={\tenrm\hfil -- \folio\ -- \hfil}
\hautpagegauche={\ifx\midfolio\oui\the\hautpagemilieu\else\tenrm\folio\hfill\the\auteurcourant\hfill\fi}
\hautpagedroite={\ifx\midfolio\oui\the\hautpagemilieu\else\hfill\the\titrecourant\hfill\tenrm\folio\fi}
\newtoks\baspagegauche \baspagegauche={\hfil}
\newtoks\baspagedroite \baspagedroite={\hfil}
\headline={\ifpagetitre\the\hautpagetitre
\else\ifodd\pageno\the\hautpagedroite\else\the\hautpagegauche\fi\fi }
\footline={\ifpagetitre\the\baspagetitre
\else\ifodd\pageno\the\baspagedroite
\else\the\baspagegauche\fi\fi \global\pagetitrefalse}

\def\pageblanche{\vfill\eject\pagetitretrue
\null\vfill\eject
\pagetitretrue
}
\def\chgtpage{\ifodd\pageno \else
\pageblanche \fi \pagetitretrue\titreun=0\footnotenumber=0}

\def\chgtpageincrtitreun{\ifodd\pageno \else
\pageblanche \fi \pagetitretrue\footnotenumber=0}

\def\majnombres{\ifodd\pageno \else
\pageblanche \fi \pagetitretrue\hautpoly\titreun=0\footnotenumber=0}

\def\hautspages#1#2{\auteurcourant={\ninepcap#1}\titrecourant={\nineit#2}}

\ifnum\chstart=\pageno \pagetitretrue\fi
  

  \def\PAR{\par}
  
  \def\leftnote#1{\vadjust{\setbox1=\vtop{\hsize 20mm\parindent=0pt\eightpoint
  \baselineskip=9pt\rightskip=4mm plus 4mm\vskip-4.7mm#1}\hbox{\kern-2cm\smash{\box1}}}}

  
  \def\raggedcenter{\leftskip=20pt plus 10em  
       \rightskip=\leftskip 
        \parfillskip=0pt 
         \spaceskip=.3333em \xspaceskip=.5em 
          \pretolerance=9999 \tolerance=9999
           \hyphenpenalty=9999 \exhyphenpenalty=9999 }
           
  \def\titrecentre#1{{\parindent0mm\raggedcenter
       \spaceskip=.6em plus .2em minus .2em\xspaceskip=.6em plus .2em minus .2em
        \tit#1\par}}
        


  \def\oui{oui}
  
\def\fontetitreun{\ifx\paradouze\oui\douzepts\gpdouze\twelvebf\textfont1=\twelveib\else
\quatorzepts\gpquatorze\fourteenbf\fi}

\def\fontetitreunl{\douzepts\textfont1=\twelveib\scriptfont1=\tenib\fourteenti}
 
 \def\fontetitredeux{\textfont1=\eleveni\ifx\paradouze\oui\onzepts\scriptfont1=\ninei\elevenit\else
                        \douzepts\twelveit\fi}
 
   \def\fontetitredeuxb{\ifx\paradouze\oui\onzepts\eleventi\gponze\textfont1=\elevenib\scriptfont1=\nineib
                         \else\douzepts\twelveti\scriptfont1=\twelveib\scriptfont1=\tenib\gpdouze\fi}
                         
\def\fontetitredeuxl{\onzepts\textfont1=\elevenbf\scriptfont1=\ninebf\twelvebf}
  
\def\fontetitretrois{\textfont0=\elevenrm\scriptfont0=\eightrm\textfont1=\eleveni
                      \scriptfont1=\eighti\scriptscriptfont1=\sixi\elevenit}
                      
\def\fontetitrequatre{\textfont0=\elevenrm\scriptfont0=\eightrm\textfont1=\eleveni
                      \scriptfont1=\eighti\scriptscriptfont1=\sixi\elevenrm}
  
  \newcount\titreun\titreun=0
  \newcount\titredeux\titredeux=0
  \newcount\titretrois\titretrois=0
  \newcount\titrequatre\titrequatre=0
  \newcount\enonce\enonce=0
  
  \def\incr#1{\global\advance#1 by 1 {\the #1}}
  \def\avance#1{\global\advance#1 by 1}
  \def\init#1{\global#1=0}
  
  \long\def\Indentation#1#2{\setbox10=\hbox{\fontetitreun#1}
  	                    \ifdim\wd10 < 4mm
                         \setbox10=\hbox to 4mm{\box10\hfill}
                       \else\ifdim\wd10 < 6mm
                         \setbox10=\hbox to 6mm{\box10\hfill}
  	                    \else\ifdim\wd10 < 8mm
                         \setbox10=\hbox to 8mm{\box10\hfill}
                       \else\ifdim\wd10 < 12mm
                         \setbox10=\hbox to 12mm{\box10\hfill}
                       \fi\fi\fi\fi
                       \dimen10=\hsize
                       \advance \dimen10 by -\wd10
                       \noindent \box10 %
                       \ignorespaces
                       \hbox{\vtop{\hsize=\dimen10\raggedright\noindent\fontetitreun#2}}}

  \long\def\paraun#1{\removelastskip\par\medskip\goodbreak\vskip0pt plus.01\vsize\penalty-100
                \vskip0pt plus-.01\vsize
  	              \init{\titredeux}\ifnum\optionparag=1{\init\eqnumber\init\enonce}\else{}\fi
                  \goodbreak{\fontetitreun
  	                \Indentation{\incr{\titreun}.\ }{\fontetitreun #1\par}}\nobreak\medskip}

 %
 %
 \long\def\paraunc#1{\removelastskip\par\bigskip\goodbreak\vskip0pt plus.01\vsize\penalty-100
                \vskip0pt plus-.01\vsize
  	              \init{\titredeux}
                 \ifnum\optionparag=1{\init{\eqnumber}\init\enonce}\else{}\fi
                  \goodbreak
  	                {\parindent0mm\raggedcenter\fontetitreun\incr{\titreun}.\ 
                     \fontetitreun #1\par}\nobreak\medskip}
                     
\newtoks\titreunl
\titreunl={\ifnum\titreun=1{I}\fi%
\ifnum\titreun=2{II}\fi%
\ifnum\titreun=3{III}\fi%
\ifnum\titreun=4{IV}\fi%
\ifnum\titreun=5{V}\fi%
\ifnum\titreun=6{VI}\fi%
\ifnum\titreun=7{VII}\fi%
\ifnum\titreun=8{VIII}\fi%
\ifnum\titreun=9{IX}\fi%
\ifnum\titreun=10{X}\fi%
\ifnum\titreun=11{XI}\fi%
\ifnum\titreun=12{XII}\fi%
\ifnum\titreun=13{XIII}\fi%
}
\long\def\paraunl#1{\removelastskip\par\bigskip\bigskip\goodbreak\vskip0pt plus.01\vsize\penalty-100
                \vskip0pt plus-.01\vsize
  	              \init{\titredeux}\ifnum\optionparag=1{\init\eqnumber\init\enonce}\else{}\fi
                  \goodbreak{\fontetitreunl
  	                \Indentation{\global\advance\titreun by 1{\the\titreunl}.\ }{\fontetitreunl #1\par}}\nobreak\smallskip}

  
  \long\def\paradeux#1{\init{\titretrois}\vskip0pt plus.01\vsize\penalty-10
                \vskip0pt plus-.01\vsize\ifx \elie\oui\medskip\ifnum\titredeux>0\medskip\fi\fi
                 \Indentation{\fontetitredeux\the\titreun${\cdot}$\incr{\titredeux}.}
                              {\fontetitredeux\textfont1=\eleveni#1}\nobreak\par }
  
  \long\def\paradeuxb#1{\init{\titretrois}\vskip0pt plus.001\vsize\penalty-10
                \vskip0pt plus-.01\vsize{\ifx \elie\oui\medskip\ifnum\titredeux>0\medskip\fi\fi
                  \Indentation
  {\fontetitredeuxb\the\titreun${\cdot}$\incr{\titredeux}.}{ \fontetitredeuxb#1}}\nobreak
\smallskip}

\newtoks\titredeuxl
\titredeuxl={\ifnum\titredeux=1{A}\fi%
\ifnum\titredeux=2{B}\fi%
\ifnum\titredeux=3{C}\fi%
\ifnum\titredeux=4{D}\fi%
\ifnum\titredeux=5{E}\fi%
\ifnum\titredeux=6{F}\fi%
\ifnum\titredeux=7{G}\fi%
\ifnum\titredeux=8{H}\fi%
\ifnum\titredeux=9{I}\fi%
\ifnum\titredeux=10{J}\fi%
\ifnum\titredeux=11{K}\fi%
\ifnum\titredeux=12{L}\fi%
\ifnum\titredeux=13{M}\fi%
}
 \long\def\paradeuxl#1{\init{\titretrois}\vskip0pt plus.001\vsize\penalty-10
                \vskip0pt plus-.01
                \vsize \bigskip%
                  \Indentation
     {\fontetitredeuxl\global\advance\titredeux by 1
  \quad \the\titreunl${\cdot}$\the\titredeuxl.}{ \fontetitredeuxl#1}
  \removelastskip\nobreak\smallskip}
  

  \long\def\paratrois#1{\init{\titrequatre}\ifdim\lastskip<\smallskipamount
                \removelastskip\smallskip\fi
                 \vskip0pt plus.01\vsize\penalty-10
                  \vskip0pt
plus-.01\vsize{\ifx \elie\oui\ifnum\titretrois>0\medskip\fi\fi
\Indentation{\fontetitretrois\the\titreun${\cdot}$\the\titredeux${\cdot}$\incr{\titretrois}.\ }
  {\hskip0mm\baselineskip=14pt\fontetitretrois#1}\nobreak\smallskip}}
  
  
  \long\def\paratroisl#1{\init{\titrequatre}\ifdim\lastskip<\smallskipamount
                \removelastskip\fi
                 \vskip0pt plus.01\vsize\penalty-10
                  \vskip0pt
plus-.01\vsize\ifx \elie\oui\bigskip
\fi
\Indentation{\fontetitretrois\quad \quad \the\titreunl{${\cdot}$}\the\titredeuxl${\cdot}$\incr{\titretrois}.\ }
  {\hskip0mm\fontetitretrois#1}\nobreak\smallskip}


  \long\def\paraquatre#1{\ifdim\lastskip<\smallskipamount
                \removelastskip\smallskip\fi
                 \vskip0pt plus.01\vsize\penalty-10
                  \vskip0pt
                  plus-.01\vsize\par
 
\Indentation{\fontetitrequatre \the\titreun{${\cdot}$}\the\titredeux${\cdot}$\the\titretrois${\cdot}$\incr{\titrequatre}.\ }
{\hskip0mm\fontetitrequatre#1}\nobreak\smallskip}


\newtoks\titrequatrel
\titrequatrel={\ifnum\titrequatre=1{a}\fi%
\ifnum\titrequatre=2{b}\fi%
\ifnum\titrequatre=3{c}\fi%
\ifnum\titrequatre=4{d}\fi%
\ifnum\titrequatre=5{e}\fi%
\ifnum\titrequatre=6{f}\fi%
\ifnum\titrequatre=7{g}\fi%
\ifnum\titrequatre=8{h}\fi%
\ifnum\titrequatre=9{i}\fi%
\ifnum\titrequatre=10{j}\fi%
\ifnum\titrequatre=11{k}\fi%
\ifnum\titrequatre=12{l}\fi%
\ifnum\titrequatre=13{m}\fi%
}
\long\def\paraquatrel#1{\ifdim\lastskip<\smallskipamount
                \removelastskip\smallskip\fi
                 \vskip0pt plus.01\vsize\penalty-10
                  \vskip0pt
                  plus-.01\vsize{\bigskip
\Indentation{\global\advance\titrequatre by 1
\fontetitrequatre\quad \quad \quad \the\titreunl${\cdot}$\the\titredeuxl${\cdot}$\the\titretrois${\cdot}$\the\titrequatrel.\ }
{\hskip0mm\fontetitrequatre#1}\nobreak\smallskip}}

\ifx\optionkeys\oui
\def\drefun#1{\definexref{¤#1}{{\the\titreun}}{}} 
\def\drefdeux#1{\definexref{¤#1}{{\the\titreun}.{\the\titredeux}}{}}
\def\dreftrois#1{\definexref{¤#1}{{\the\titreun}.{\the\titredeux}.{\the\titretrois}}{}}
\else
\def\drefun#1{\definexref{prg#1}{{\the\titreun}}{}} 
\def\drefdeux#1{\definexref{prg#1}{{\the\titreun}.{\the\titredeux}}{}}
\def\dreftrois#1{\definexref{prg#1}{{\the\titreun}.{\the\titredeux}.{\the\titretrois}}{}}
\fi

%


  \long\def\propdeux#1#2#3#4{%
       \avance{\enonce}
       \leavevmode\edef\temp{#2}%
         \ifx\temp\empty 
          \else
           \definexref{#2}{#1~{\the\titreun.\the\enonce}}{enonces}
            \definexref{s#2}{{\the\titreun.\the\enonce}}{enonces}
             \fi
\smallskip
      \noindent{\bf#1\ {\bf\the\titreun.\the\enonce{#3}.}\enspace}{\sl#4\par}%
      \ifdim\lastskip<\medskipamount \removelastskip\penalty55\par \fi
   }

  \long\def\propun#1#2#3#4{%
      \avance{\enonce}
       \leavevmode\edef\temp{#2}%
        \ifx\temp\empty 
          \else
           \definexref{#2}{#1~{\the\enonce}}{enonces}
            \definexref{{s#2}}{{\the\enonce}}{enonces}
             \fi
   \par 
     \noindent{\bf#1\ {\bf\the\enonce{#3}.}\enspace}{\sl#4\par}%
     \ifdim\lastskip<\medskipamount \removelastskip\penalty55\medskip\fi
  }
  
  \long\def\prop#1#2#3#4{\ifnum\optionparag=1
                          \propdeux{#1}{#2}{\textfont1=\elevenib#3}{#4} \else\propun{#1}{#2}{\textfont1=\elevenib#3}{#4}\fi}

  \long\def\propt#1#2#3{\ifx\tpf\oui \prop{Th\'eo\-r\`eme}{#1}{#2}{#3}\par
                       \else\prop{Theorem}{#1}{#2}{#3}\par\fi}
  \long\def\Propt#1#2{\propt{#1}{}{#2}}
  \long\def\propl#1#2#3{\ifx\tpf\oui\prop{Lem\-me}{#1}{#2}{#3}\par
                         \else\prop{Lemma}{#1}{#2}{#3}\par\fi}
  \long\def\Propl#1#2{\propl{#1}{}{#2}}
  \long\def\propc#1#2#3{\ifx\tpf\oui\prop{Corol\-laire}{#1}{#2}{#3}\par
                         \else\prop{Corollary}{#1}{#2}{#3}\par\fi}
  
  \long\def\propp#1#2#3{\prop{Pro\-po\-si\-tion}{#1}{#2}{#3}\par}
  \long\def\Propp#1#2{\propp{#1}{}{#2}} 
  \long\def\propd#1#2#3{\ifx\tpf\oui\prop{D\'efi\-nition}{#1}{#2}{#3}\par
                       \else\prop{Definition}{#1}{#2}{#3}\par\fi} 
  
  \long\def\proptd#1#2#3{\ifx\tpf\oui\prop{Th\'eor\`eme et d\'efi\-nition}{#1}{#2}{#3}\par
                       \else\prop{Theorem and definition}{#1}{#2}{#3}\par\fi}


  
  \newcount\optionparag\optionparag=1
  
  \long\def\section#1#2{\ifnum\optionparag=1 \paraun{#2} 
                        \else\goodbreak{\fontetitreun
  	                \Indentation{#1.\ }{#2}}\nobreak\smallskip\fi}

  \def\eqconstruct#1{\ifnum\optionparag=1{\the\titreun\hbox{$\cdot$}#1}\else{#1}\fi}

  
  
  \def\numref{oui}  
  
  \newcount\mesref\mesref=0 
  \def\defbib#1{\ifx\numref\oui\global\advance\mesref by 1 \definexref{#1}{{\the
                 \mesref}}{}\else\definexref{#1}{#1}{}\fi}
  \def\bibtem#1{\defbib{#1}\item{\citer{#1}}}
  \def\citer#1{[\ref{#1}]}
  \def\citeplus#1#2{[\ref{#1}; #2]}

  
  \font\seventeenmsa=msam10 at 17pt    
  \font\fourteenmsa=msam10 at 14pt
  \font\twelvemsa=msam10 at 12pt
  \font\tenmsa=msam10                 
  \font\ninemsa=msam10 at 9pt 
  \font\eightmsa=msam10 at 8pt 
  \font\sevenmsa=msam7 
  \font\sixmsa=msam10 at 6pt
  \font\fivemsa=msam5
  \newfam\msafam\textfont\msafam=\tenmsa\scriptfont\msafam=\sevenmsa\scriptscriptfont\msafam=\fivemsa
  
  \font\seventeenbb=msbm10 at 17pt     
  \font\fourteenbb=msbm10 at 14pt
  \font\twelvebb=msbm10 at 12pt
  \font\tenbb=msbm10                   
  \font\ninebb=msbm10 at 9pt 
  \font\eightbb=msbm10 at 8pt 
  \font\sevenbb=msbm7 
  \font\sixbb=msbm10 at 6pt
  \font\fivebb=msbm5 
  \newfam\bbfam\textfont\bbfam=\tenbb\scriptfont\bbfam=\sevenbb\scriptscriptfont\bbfam=\fivebb
  \def\bb{\fam\bbfam\tenbb}%

  \font\seventeenscaln=eusm10 at 17pt   
  \font\twelvescaln=eusm10 at 12pt
  \font\tenscaln=eusm10                
  \font\ninescaln=eusm10 scaled 900
  \font\eightscaln=eusm10 scaled 800
  \font\sevenscaln=eusm10 scaled 700
  \font\sixscaln=eusm10 scaled 600
   
  \newfam\scalnfam\textfont\scalnfam=\tenscaln\scriptfont\scalnfam=\sevenscaln\scriptscriptfont\scalnfam=\sixscaln
  \def\scaln{\fam\scalnfam\tenscaln}%
  \def\scal{\scaln}
  
  \font\tenscalb=eusb10                

  \font\sevenscalb=eusb10 scaled 700

  \newfam\scalbfam\textfont\scalbfam=\tenscalb\scriptfont\scalbfam=\sevenscalb
  %
  
  %
  %
  \font\fourteenrm=cmr12 scaled 1200
  \font\elevenrm=cmr10 at 11pt
  \font\twelverm=cmr12
  \font\ninerm=cmr9
  \font\eightrm=cmr8      
  \font\sevenrm=cmr7
  \font\sixrm=cmr6

  \font\seventeenpcap=cmcsc10 at 17pt
  \font\tenpcap=cmcsc10                        
  \font\ninepcap=cmcsc9
  \font\eightpcap=cmcsc8
  \font\sevenpcap=cmcsc10 scaled 700
  
  \newfam\pcapfam\textfont\pcapfam=\tenpcap\scriptfont\pcapfam=\sevenpcap
  \def\pcap{\fam\pcapfam\tenpcap}
  
  \font\seventeenrm=cmbx12 scaled 1400

  \font\fourteenbf=cmbx10 scaled 1400
  
  \font\twelvebf=cmbx12
  \font\elevenbf=cmbx10 at 11pt
  \font\ninebf=cmbx9  
  \font\eightbf=cmbx8
  \font\sixbf=cmbx6
  
  \font\tengot=eufm10                           
   
  \font\eightgot=eufm10 at 8truept 
  \font\sevengot=eufm7 
  \font\sixgot=eufm10 at 6 truept 
   
  \newfam\gotfam
  \textfont\gotfam=\tengot\scriptfont\gotfam=\sevengot\scriptscriptfont\gotfam=\sixgot
  \def\got{\fam\gotfam\tengot}%

  
  \def\tit{%
  \textfont0=\seventeenrm\scriptfont0=\tenrm\def\rm{\fam0\seventeenrm}%
  \textfont1=\seventeenib\scriptfont1=\twelveib%
  \textfont2=\seventeensy\scriptfont2=\twelvesy\scriptscriptfont2=\ninesy
  \textfont3=\seventeenex
  \textfont\itfam=\seventeenti
  \def\it{\fam\itfam\seventeenti}%
  \textfont\bbfam=\seventeenbb \scriptfont\bbfam=\twelvebb
  \def\bb{\fam\bbfam\seventeenbb}%
  \textfont\msafam=\seventeenmsa\scriptfont\msafam=\twelvemsa
  \textfont\scalnfam=\seventeenscaln
  \def\pcap{\fam\pcapfam\seventeenpcap}
  \normalbaselineskip=25pt\normalbaselines\rm}

  \font\seventeenti=cmbxti10 scaled 1680
  
  \font\fourteenti=cmbxti10 at 14pt
  
  \font\twelveti=cmbxti10 scaled 1200
  \font\eleventi=cmbxti10 at 11pt

  %
  %
  \font\twelveit=cmti12	
  \font\elevenit=cmti10 scaled 1100
  \font\nineit=cmti9
  \font\eightit=cmti8
  \font\sevenit=cmti7

  %
  %
 
 \font\seventeenib=cmmib10 scaled 1680
  \font\fourteenib=cmmib10 scaled 1400
  \font\twelveib=cmmib10 scaled 1200
  \font\elevenib=cmmib10 scaled 1100
  \font\tenib=cmmib10
\font\eightib=cmmib10 scaled 800
  \font\nineib=cmmib10 scaled 900
\font\sevenib=cmmib10 scaled 700
\font\sixib=cmmib10 scaled 600
\font\fiveib=cmmib10 scaled 500

\ifx\ITAN\oui
\else
\innernewfam\cmmibfam
\textfont\cmmibfam=\tenib
\scriptfont\cmmibfam=\sevenib
\scriptscriptfont\cmmibfam=\fiveib
\def\ib{\fam\cmmibfam\tenib}
\fi

  %
  %
  \font\twelvei=cmmi10 scaled 1200
  \font\eleveni=cmmi10 scaled 1100
  \font\ninei=cmmi9
  \font\eighti=cmmi8 
  \font\seveni=cmmi7 			                
  \font\sixi=cmmi6
  
  \font\ninesl=cmsl9                    
  \font\eightsl=cmsl8 
  \font\sevensl=cmsl10 at 7pt

  \font\ninett=cmtt9                    
  \font\eighttt=cmtt8
  \font\seventt=cmtt10 scaled 700

  \font\seventeensy=cmsy10 scaled 1680    
  \font\fourteensy=cmsy10 scaled 1400
  \font\twelvesy=cmsy10 scaled 1176
  
  \font\ninesy=cmsy9                      
  \font\eightsy=cmsy8
  \font\sixsy=cmsy6
  \font\seventeenex=cmex10 at 17pt
  \font\fourteenex=cmex10 at 14pt
  \font\twelveex=cmex10 at 12pt
  \font\nineex=cmex10 at 9pt
  \font\eightex=cmex10 at 8pt
  \font\sevenex=cmex10 at 7pt
  \font\sixex=cmex10 at 6pt
  \font\fiveex=cmex10 at 5pt
  
   
  \font\fourteengp=cmmi10 at 14pt
  
  \font\twelvegp=cmmib10 at 12pt
  \font\elevengp=cmmib10 at 11pt
  \font\tengp=cmmib10                          
  \font\ninegp=cmmib10 at 9pt 
  \font\eightgp=cmmib8 
  \font\sevengp=cmmib7 
  \font\sixgp=cmmib6


  \def\gponze{\textfont0=\elevenbf\scriptfont0=\eightbf\scriptscriptfont0=\sixbf
           \textfont1=\elevengp\scriptfont1=\eightgp\scriptscriptfont1=\sixgp}
  \def\gpdouze{\textfont0=\twelvebf\scriptfont0=\tenbf\scriptscriptfont0=\ninebf
           \textfont1=\twelvegp\scriptfont1=\tengp\scriptscriptfont1=\ninegp}        
  
 \def\gpquatorze{\textfont0=\fourteenbf\scriptfont0=\twelvebf\scriptscriptfont0=\elevenbf
           \textfont1=\fourteengp\scriptfont1=\twelvegp\scriptscriptfont1=\elevengp}

  
  \expandafter\chardef\csname pre amssym.def at\endcsname=\the\catcode`\@
  \catcode`\@=11
  \def\undefine#1{\let#1\undefined}
  \def\newsymbol#1#2#3#4#5{\let\next@\relax
   \ifnum#2=\@ne\let\next@\msafam@\else
   \ifnum#2=\tw@\let\next@\bbfam@\fi\fi
   \mathchardef#1="#3\next@#4#5}
  \def\mathhexbox@#1#2#3{\relax
   \ifmmode\mathpalette{}{\m@th\mathchar"#1#2#3}%
   \else\leavevmode\hbox{$\m@th\mathchar"#1#2#3$}\fi}
  \def\hexnumber@#1{\ifcase#1 0\or 1\or 2\or 3\or 4\or 5\or 6\or 7\or 8\or
   9\or A\or B\or C\or D\or E\or F\fi}
  
  \def\setboxz@h{\setbox\z@\hbox}
  \def\wdz@{\wd\z@}
  \def\boxz@{\box\z@}
  
  \edef\msafam@{\hexnumber@\msafam}
  \mathchardef\dabar@"0\msafam@39
  
  \edef\bbfam@{\hexnumber@\bbfam}
  \def\widehat#1{\setboxz@h{$\m@th#1$}%
   \ifdim\wdz@>\tw@ em\mathaccent"0\bbfam@5B{#1}%
   \else\mathaccent"0362{#1}\fi}
  \def\widetilde#1{\setboxz@h{$\m@th#1$}%
   \ifdim\wdz@>\tw@ em\mathaccent"0\bbfam@5D{#1}%
   \else\mathaccent"0365{#1}\fi}
  \newsymbol\leqq 1335          
  \newsymbol\leqslant 1336
  \newsymbol\lessgtr 1337       
  \newsymbol\backprime 1038     
  \newsymbol\risingdotseq 133A  
  \newsymbol\fallingdotseq 133B 
  \newsymbol\succcurlyeq 133C   
  \newsymbol\geqq 133D          
  \newsymbol\geqslant 133E
  \newsymbol\nmid 232D
  \newsymbol\nexists 2040
  \newsymbol\smallsetminus 2272
  \newsymbol\varnothing 203F 
  \catcode`\@=\active

  \catcode`\@=11
  \newcount\typofr\typofr=1
  
  \catcode`\;=\active
  \def;{\ifnum\typofr=1\relax\ifhmode\ifdim\lastskip>\z@\unskip\fi
     \kern.2em\fi\string;\else\string;\fi}
  
  \catcode`\:=\active
  \def:{\ifnum\typofr=1\relax\ifhmode\ifdim\lastskip>\z@\unskip\fi
  \penalty\@M\ \fi\string:\else\string:\fi}
  
  \catcode`\!=\active
  \def!{\ifnum\typofr=1\relax\ifhmode\ifdim\lastskip>\z@\unskip\fi
     \kern.2em\fi\string!\else\string!\fi}
  
  \catcode`\?=\active
  \def?{\ifnum\typofr=1\relax\ifhmode\ifdim\lastskip>\z@\unskip\fi
     \kern.2em\fi\string?\else\string?\fi}

  \def\francais{\typofr=1\def\tpf{oui}}
  \def\anglais{\typofr=2\def\tpf{non}\def\english{oui}}
  \def\oui{oui}
  \francais
  
  \catcode`\@=12
  

\ifx\textures\oui
\def\raye #1|{\leavevmode\setbox1=\hbox{#1}%
\raise .5pt\hbox to \wd1{\xleaders\hbox{\rge{$ \sct / $}%
\kern 1pt}\hfill\kern -1pt }\kern -\wd1 \unhbox1\relax }

\def\barre #1|{\leavevmode\setbox1=\hbox{#1}%
\rlap{\Red\vrule height 2.4pt depth -1.2pt width \wd1}\Black \unhbox1\relax}
\else
\def\raye #1|{\leavevmode\setbox1=\hbox{#1}%
\raise .5pt\hbox to \wd1{\xleaders\hbox{\rge{$ \sct / $}%
\kern 1pt}\hfill\kern -1pt }\kern -\wd1 \unhbox1\relax }

\def\barre #1|{\leavevmode\setbox1=\hbox{#1}%
\rlap{\color{red}\vrule height 2.4pt depth -1.2pt width \wd1}\color{black} \unhbox1\relax}

\fi
  

  
  \def\og{\leavevmode\raise.24ex\hbox{$\scriptscriptstyle\langle\!\langle\>$}}    
  \def\fg{\leavevmode\raise.24ex\hbox{$\scriptscriptstyle\>\rangle\!\rangle$}}    

  \def\d{\,{\rm d}}
  \def\dd{{\rm d}}

  \def\z{{\bb Z}}
  \def\r{{\bb R}}
  
  \def\N{{\bb N}}

  \def\PP{{\bb P}}

  \def\HH{{\scal H}}

  \def\O{{\scal O}}
  \def\P{{\scaln P}}

  \def\Y{{\scal Y}}

  \def\frac#1#2{{#1\over #2}}
  \def\di#1#2{\sct#1\atop{\sct#2}}

  \def\qedbox{$\rlap{$\sqcap$}\sqcup$}           
  \def\qed{\nobreak\hfill\penalty250 \hbox{}\nobreak\hfill\qedbox\par }

  \def\¤{\S\thinspace}

  \def\¥{$\bullet$ }
  
  
  \def\e{{\rm e}}
  
  \def\md#1#2{\equiv#1\,({\rm mod\,}#2)}

  \def\epsilon{\varepsilon}

  \def\phi{\varphi}
  \def\theta{\vartheta}
  \def\rho{\varrho}
  \def\dm{{\textstyle{1\over 2}}}

  \def\sct{\scriptstyle}
  \def\pf{\noi{\it Proof. }}
  \def\nid{\ifnum\typofr=1\par\noindent{\it D\'emonstration. }\else\pf\fi}
  \def\noi{\noindent}
  \def\rem{\ifnum\typofr=1\noi{\it Remarque.}\ \else\noi{\it Remark.}\ \fi}
  \def\rems{\ifnum\typofr=1\noi{\it Remarques.}\ \else\noi{\it Remarks.}\ \fi}

  \def\sset{\smallsetminus}

  \def\1{{\bf 1}}
  \def\|{\Vert}

  \def\leq{\leqslant}
  \def\ge{\geqslant}\def\geq{\geqslant}

  \def\eg{{e.g.}}
  \newsymbol\subsetneqq 2324
  \newsymbol\subsetneq 2328

  \def\fl#1{\left\lfloor #1 \right\rfloor}

  \def\log{\mathop{\rm log}\nolimits}




  \def\pmb#1{\setbox0=\hbox{#1}%
  \kern-.025em\copy0\kern-\wd0\kern.05em\copy0\kern-\wd0\kern-.025em\raise .0433em\box0 }

  
  \skewchar\eighti='177 \skewchar\sixi='177
  \skewchar\eightsy='60 \skewchar\sixsy='60
  
  \def\eightpoint{%
  \textfont0=\eightrm\scriptfont0=\sixrm\scriptscriptfont0=\fiverm
  \def\rm{\fam0\eightrm}%
  \textfont1=\eighti\scriptfont1=\sixi
  \scriptscriptfont1=\fivei\def\oldstyle{\fam1\seveni}%
  \textfont2=\eightsy\scriptfont2=\sixsy\scriptscriptfont2=\fivesy
  \textfont3=\eightex\scriptfont3=\sixex
  \textfont\itfam=\eightit
  \def\it{\fam\itfam\eightit}%
  \textfont\slfam=\eightsl
  \def\sl{\fam\slfam\eightsl}%
  \textfont\bbfam=\eightbb \scriptfont\bbfam=\sixbb\scriptscriptfont\bbfam=\fivebb
  \def\bb{\fam\bbfam\eightbb}%
  \textfont\msafam=\eightmsa\scriptfont\msafam=\sixmsa
  \textfont\scalnfam=\eightscaln
  \def\scaln{\fam\scalnfam\eightscaln}
  \textfont\ttfam=\eighttt
  \def\tt{\fam\ttfam\eighttt}%
\textfont\gotfam=\eightgot
  \textfont\bffam=\eightbf\scriptfont\bffam=\sixbf\scriptscriptfont\bffam=\fivebf
  \def\bf{\fam\bffam\eightbf}%
  \ifx\ITAN\oui\else\textfont\cmmibfam=\eightib
       \scriptfont\cmmibfam=\sixib
        \scriptscriptfont\cmmibfam=\fiveib
         \def\ib{\fam\cmmibfam\eightib}
   \fi
  \textfont\pcapfam=\eightpcap
  \def\pcap{\fam\pcapfam\eightpcap}
  \abovedisplayskip=2pt plus2pt minus 2pt
  \belowdisplayskip=2pt plus1pt minus 2pt
  \abovedisplayshortskip= 1pt plus 2pt minus 1pt
  \belowdisplayshortskip= 1pt plus 2pt minus 1pt
  \smallskipamount=2pt plus 1pt minus 2pt
  \medskipamount=3pt plus 2pt minus 2pt
  \bigskipamount=7pt plus 3pt minus 3pt
  \setbox\strutbox=\hbox{\vrule height 5pt depth 2pt width 0pt}%
  \normalbaselineskip=9pt\normalbaselines\rm}

  \def\({\left(}
  \def\){\right)}
  
  \def\footnoterule{\kern -2pt\hrule width 7truecm\kern 2.4pt}
  
  \def\xnotedef#1{\definexref{#1}{\noexpand\number\footnotenumber}{Note}}%

  
  
  \def\ninepoint{%
  \textfont0=\ninerm\scriptfont0=\sixrm\scriptscriptfont0=\fiverm
  \def\rm{\fam0\ninerm}%
  \textfont1=\ninei\scriptfont1=\sixi
  \scriptscriptfont1=\fivei\def\oldstyle{\fam1\ninei}%
  \textfont2=\ninesy\scriptfont2=\sixsy\scriptscriptfont2=\fivesy
  \textfont3=\nineex\scriptfont3=\sixex
  \textfont\itfam=\nineit
  \def\it{\fam\itfam\nineit}%
  \textfont\slfam=\ninesl
  \def\sl{\fam\slfam\ninesl}%
  \textfont\bbfam=\ninebb\scriptfont\bbfam=\sixbb\scriptscriptfont\bbfam=\fivebb
  \def\bb{\fam\bbfam\ninebb}%
  \textfont\msafam=\ninemsa\scriptfont\msafam=\sixmsa\scriptscriptfont\msafam=\fivemsa
  \textfont\scalnfam=\ninescaln
  \def\scaln{\fam\scalnfam\ninescaln}
  \textfont\ttfam=\ninett
  \def\tt{\fam\ttfam\ninett}%
  \textfont\bffam=\ninebf\scriptfont\bffam=\sixbf\scriptscriptfont\bffam=\fivebf
  \def\bf{\fam\bffam\ninebf}%
  \abovedisplayskip=3pt plus2pt minus 2pt
  \belowdisplayskip=3pt plus1pt minus 2pt
  \abovedisplayshortskip= 2pt plus 2pt minus 1pt
  \belowdisplayshortskip= 2pt plus 2pt minus 1pt
  \smallskipamount=2pt plus 1pt minus 2pt
  \medskipamount=3pt plus 2pt minus 2pt
  \bigskipamount=7pt plus 3pt minus 3pt
  \setbox\strutbox=\hbox{\vrule height 5pt depth 2pt width 0pt}%
  \normalbaselineskip=11pt plus.3pt minus.2pt\normalbaselines\rm}

  \def\sevenpoint{%
  \textfont0=\sevenrm\scriptfont0=\sixrm\scriptscriptfont0=\fiverm
  \def\rm{\fam0\sevenrm}%
  \textfont1=\seveni\scriptfont1=\sixi
  \scriptscriptfont1=\fivei\def\oldstyle{\fam1\seveni}%
  \textfont2=\sevensy\scriptfont2=\sixsy\scriptscriptfont2=\fivesy
  \textfont3=\sevenex\scriptfont3=\fiveex
  \textfont\itfam=\sevenit
  \def\it{\fam\itfam\sevenit}%
  \textfont\slfam=\sevensl
  \def\sl{\fam\slfam\sevensl}%
  \textfont\bbfam=\sevenbb \scriptfont\bbfam=\sixbb\scriptscriptfont\bbfam=\fivebb
  \def\bb{\fam\bbfam\sevenbb}%
  \textfont\msafam=\sevenmsa\scriptfont\msafam=\sixmsa
  \textfont\scalnfam=\sevenscaln
  \def\scaln{\fam\scalnfam\sevenscaln}
  \textfont\bffam=\sevenbf\scriptfont\bffam=\sixbf\scriptscriptfont\bffam=\fivebf
  \def\bf{\fam\bffam\sevenbf}%
  \textfont\ttfam=\seventt
  \abovedisplayskip=2pt plus2pt minus 2pt
  \belowdisplayskip=2pt plus1pt minus 2pt
  \abovedisplayshortskip= 1pt plus 2pt minus 1pt
  \belowdisplayshortskip= 1pt plus 2pt minus 1pt
  \smallskipamount=2pt plus 1pt minus 2pt
  \medskipamount=3pt plus 2pt minus 2pt
  \bigskipamount=7pt plus 3pt minus 3pt
  \setbox\strutbox=\hbox{\vrule height 5pt depth 2pt width 0pt}%
  \normalbaselineskip=9pt\normalbaselines\rm}

 \def\onzepts{%
 \textfont0=\elevenrm\scriptfont0=\ninerm
 \textfont1=\eleveni\scriptfont1=\ninei
}

\def\douzepts{%
  \textfont0=\twelverm\scriptfont0=\tenrm\def\rm{\fam0\twelverm}%
  \textfont1=\twelvei\scriptfont1=\teni%
  \textfont2=\twelvesy\scriptfont2=\tensy\scriptscriptfont2=\eightsy
  \textfont3=\twelveex
  \textfont\itfam=\twelveti
  \def\it{\fam\itfam\twelveti}%
  \textfont\bffam=\twelvebf\scriptfont\bffam=\tenbf\scriptscriptfont\bffam=\eightbf
  \def\bf{\fam\bffam\twelvebf}%
  \textfont\bbfam=\twelvebb \scriptfont\bbfam=\tenbb
  \def\bb{\fam\bbfam\twelvebb}%
  \textfont\msafam=\twelvemsa\scriptfont\msafam=\tenmsa
  \textfont\scalnfam=\twelvescaln
  \normalbaselineskip=15pt\normalbaselines\rm}

\def\quatorzepts{%
  \textfont0=\fourteenrm\scriptfont0=\twelverm\def\rm{\fam0\fourteenrm}%
  \textfont1=\fourteenib\scriptfont1=\twelveib%
  \textfont2=\fourteensy\scriptfont2=\twelvesy\scriptscriptfont2=\tensy
  \textfont3=\fourteenex
  \textfont\itfam=\fourteenti
  \def\it{\fam\itfam\fourteenti}%
  \textfont\bffam=\fourteenbf\scriptfont\bffam=\twelvebf\scriptscriptfont\bffam=\tenbf
  \def\bf{\fam\bffam\fourteenbf}%
  \textfont\bbfam=\fourteenbb \scriptfont\bbfam=\twelvebb
  \def\bb{\fam\bbfam\fourteenbb}%
  \textfont\msafam=\fourteenmsa\scriptfont\msafam=\twelvemsa
  \textfont\scalnfam=\twelvescaln
  \normalbaselineskip=18pt\normalbaselines\rm}


\def\AA{{\it Acta Arith.}}

\def\picture #1 by #2 (#3){\leavevmode\vbox to #2{
     \hrule width #1 height 0pt depth 0pt
      \vfill
       \special{picture #3}}}

\def\illustration #1 by #2 (#3) scaled #4{\dimen1=#2
  \divide\dimen1 by 1000
  \multiply\dimen1 by #4
  \vtop to \dimen1{\dimen1=#1
  \divide\dimen1 by 1000
  \multiply\dimen1 by #4
  \hsize=\dimen1\vss
  \noindent\special{illustration #3 scaled #4}}}

\ifx\couleurs\oui

\fi

\anglais
\beginpackages
\endpackages
\optionparag=1
\def\paradouze{oui}
\vsize=250truemm
\voffset=-3truemm
\ifx\optionkeymacros\undefined\else \fi

\catcode`\Œ=\active\defŒ{{\aa}}       
\catcode`\º=\active\defº{\int}        
\catcode`\=\active\def{\c c}        
\catcode`\¶=\active\def¶{\partial}    
\catcode`\Ä=\active\defÄ{\oint}       
\catcode`\Æ=\active\defÆ{\triangle}   
\catcode`\Â=\active\defÂ{\neg}        
\catcode`\µ=\active\defµ{\mu}         
\catcode`\¿=\active\def¿{{\o}}        
\catcode`\¹=\active\def¹{\pi}         
\catcode`\Ï=\active\defÏ{{\oe}}       
\catcode`\§=\active\def§{{\ss}}       
\catcode`\ =\active\def {\dagger}     
\catcode`\Ã=\active\defÃ{\sqrt}       
\catcode`\·=\active\def·{\Sigma}      
\catcode`\Å=\active\defÅ{\approx}     
\catcode`\½=\active\def½{\Omega}      
\catcode`\£=\active\def£{{\it\$}}     
\catcode`\°=\active\def°{\infty}      
\catcode`\¤=\active\def¤{{\S}}        
\catcode`\¦=\active\def¦{{\P}}        
\catcode`\¥=\active\def¥{\bullet}     
\catcode`\»=\active\def»{\leavevmode\raise.585ex\hbox{\b a}}      
\catcode`\¼=\active\def¼{\leavevmode\raise.6ex\hbox{\b o}}        
\catcode`\­=\active\def­{\not=}       
\catcode`\²=\active\def²{\leq}        
\catcode`\³=\active\def³{\geq}        
\catcode`\Ö=\active\defÖ{\div}        
\catcode`\É=\active\defÉ{{\dots}}     
\catcode`\¾=\active\def¾{{\ae}}       
\catcode`\Ç=\active\defÇ{\og}         
\catcode`\Ò=\active\defÒ{``}          
\catcode`\Á=\active\defÁ{!`}          
\catcode`\¢=\active\def¢{\rlap/c}     
\catcode`\Ô=\active\defÔ{`}           
\catcode`\Õ=\active\defÕ{'}           


\catcode`\=\active\def{{\AA}}       
\catcode`\'=\active\def'{\c C}        
\catcode`\¯=\active\def¯{{\O}}        
\catcode`\¸=\active\def¸{\Pi}         
\catcode`\Î=\active\defÎ{{\OE}}       
\catcode`\®=\active\def®{{\AE}}       
\catcode`\×=\active\def×{\diamond}    
\catcode`\¡=\active\def¡{\accent'27}  
\catcode`\Ó=\active\defÓ{''}          
\catcode`\±=\active\def±{\pm}         
\catcode`\È=\active\defÈ{\fg}         
\catcode`\À=\active\defÀ{?`}          
\catcode`\Ð=\active\defÐ{--}          
\catcode`\Ñ=\active\defÑ{---}         


\catcode`\Š=\active\defŠ{\"a}        
\catcode`\'=\active\def'{\"e}        
\catcode`\•=\active\def•{\"{\i}}     
\catcode`\š=\active\defš{\"o}        
\catcode`\Ÿ=\active\defŸ{\"u}        
\catcode`\Ø=\active\defØ{\"y}        
\catcode`\å=\active\defå{\^A}        
\catcode`\€=\active\def€{\"A}        
\catcode`\…=\active\def…{\"O}        
\catcode`\†=\active\def†{\"U}        
\catcode`\‡=\active\def‡{\'a}        
\catcode`\Ž=\active\defŽ{\'e}        
\catcode`\'=\active\def'{\'{\i}}     
\catcode`\—=\active\def—{\'o}        
\catcode`\œ=\active\defœ{\'u}        
\catcode`\ƒ=\active\defƒ{\'E}        
\catcode`\æ=\active\defæ{\^E}        
\catcode`\é=\active\defé{\`E}        %
\catcode`\ˆ=\active\defˆ{\`a}        
\catcode`\=\active\def{\`e}        
\catcode`\"=\active\def"{\`{\i}}     
\catcode`\˜=\active\def˜{\`o}        
\catcode`\=\active\def{\`u}        
\catcode`\Ë=\active\defË{\`A}        
\catcode`\‹=\active\def‹{\~a}        
\catcode`\–=\active\def–{\~n}        
\catcode`\›=\active\def›{\~o}        
\catcode`\Ì=\active\defÌ{\~A}        
\catcode`\"=\active\def"{\~N}        
\catcode`\Í=\active\defÍ{\~O}        
\catcode`\‰=\active\def‰{\^a}        
\catcode`\=\active\def{\^e}        
\catcode`\"=\active\def"{\^{\i}}     
\catcode`\™=\active\def™{\^o}        
\catcode`\ž=\active\defž{\^u}        

\let\optionkeymacros\null

\def\EE{{\bb E}}

\def\sE{{\scal E}}

\font\tengp=cmmib10

\font\ninegp=cmmib10 at 9pt
               \font\eightgp=cmmib8
              \font\sevengp=cmmib7
               
               \newfam\gpfam
               \textfont\gpfam=\tengp\scriptfont\gpfam=\sevengp

\font\tenib=cmmib10
             \font\nineib=cmmib10 scaled 900
\font\sevenib=cmmib10 scaled 700
             \font\fiveib=cmmib10 scaled 500

\def\gb{{\got b}}

\def\gd{{\got d}}

\def\gq{{\got q}}

\def\som#1{\mathop{\sum^{#1}}}
\def\auteur{R. de la Bretche \& G. Tenenbaum}
 
\def\titrart{Note on the mean value of the Erd\H os--Hooley Delta-function}

\hautspages{\auteur}{\titrart}
\dateheure
\titrecentre{\titrart}
\bigskip\medskip
\centerline {\auteur} 
\bigskip
{\leftskip80mm\obeylines 
\it To the memory of Richard R. Hall,
whose talent and elegance
will remain as a constant 
source of inspiration \par }
\bigskip
{\eightpoint\leftskip1cm\rightskip1cm
\noi{\bf Abstract.}  
For integer $n\geqslant 1$ and real $u$, let $\Delta(n,u):=|\{d:d\mid n,\,\e^u<d\leqslant \e^{u+1}\}|$. The Erd\H os--Hooley Delta-function is then defined by $\Delta(n):=\max_{u\in\r}\Delta(n,u).$ We improve a recent upper bound for the mean value of this function by showing that, for large $x$, we have $$\sum_{n\leqslant x}\Delta(n)\ll x(\log_2x)^{ 5/2}.$$\PAR
\medskip
\noi
{\bf Keywords:}  Distribution of divisors, concentration of divisors, average order, Hooley's function, Erd\H os-Hooley's Delta-function, Waring's problem.\par
\smallskip 
\noi \bf 2020 Mathematics Subject Classification: \rm primary   11N37; secondary 11K65.\par }
\bigskip
\medskip
\paraun{Introduction and statement of results} 
 For integer $n\geqslant 1$ and real $u$, put $$\Delta(n,u):=|\{d:d\mid n,\,\e^u<d\leqslant \e^{u+1}\}|,\quad\Delta(n):=\max_{u\in\r}\Delta(n,u).$$
 Introduced by Erd\H os \citer{Er73} (see also \citer{EN75}) and studied by Hooley \citer{Ho79}, the $\Delta$-function and specifically its (possibly weighted) mean-value proved very useful in several branches of number theory --- see, \eg,  \citer{Te86}, \citer{HT88}, and \citer{BT23}  for further references. If $\tau(n)$ denotes the total number of divisors of $n$, then $\Delta(n)/\tau(n)$ coincides with the concentration of the numbers $\log d$, $d|n$.
 \par 
 Asymptotic estimates for 
 $$S(x):=\sum_{n\leqslant x}\Delta(n)$$
 have a rather long history since Hooley's pioneer work \citer{Ho79}: see \citer{HT84}, \citer{Te85}, \citer{HT86}, \citer{HT88}, and the recent papers \citer{BT23}, \citer{KT23}, and \citer{FKT23} for a description of the main steps. While Hooley's upper bound was $S(x)\ll x(\log x)^{4/\pi-1}$ and following works aimed at reducing the value of the exponent, the first estimate of the type 
 $$S(x)\ll x\e^{c\sqrt{\log_2 x\log_3x}}\qquad (x\geqslant 16)\eqdef{GT85}$$
 appears in Tenenbaum's paper \citer{Te85}. Here and in the sequel, we let $\log _k$ denote the $k$-fold iterated logarithm. In \citer{BT23}, we could remove the $\log_3x$ from \eqref{GT85}. Very recently Koukoulopoulos and Tao made a further breakthrough by showing
 $$S(x)\ll x(\log_2x)^{11/4}\qquad (x\geqslant 3).$$
In this work, we slightly  improve  this upper bound. 
\par 
Considering lower bounds, it is shown in \citer{HT82} that $S(x)\gg x\log_2x$ and this remained unimproved for over four decades. A very recent work of Ford, Koukoulopoulos and Tao~\citer{FKT23} provides the lower bound $$S(x)\gg x(\log_2x)^{1+\eta+o(1)}\eqdef{FKT}$$
where $\eta\approx0.353327$ is the exponent appearing in the new lower bound for the normal order of $\Delta(n)$ by Ford, Green and Koukoulopoulos \citer{FGK23}. 
\par\goodbreak
\Propt{thD}{For all large real $x$, we have
$$ S(x)
\ll x(\log_2x)^{5/2}.\eqdef{estsD}$$}
Since the work \citer{Te85}, all upper estimates for $S(x)$ rest on properties of moments
$$M_q(n):=\int_\r\Delta(n,u)^q\d u\qquad (n\geqslant 1,\,q\geqslant 1).$$ As we clearly have  $\Delta(np,u)=\Delta(n,u)+\Delta(n,u-\log p)$ for any integer $n\geqslant 1$ and any prime number $p$ such that $p\nmid n$, we may write the inductive inequality
$$2M_q(n)\leqslant M_q(np)\leqslant2M_q(n)+2W_q(n,p)\quad(p\nmid n)\eqdef{recMq}$$
with
$$\eqalign{N_{j,q}(n,p)&:=\int_\r\Delta(n,u)^j\Delta(n,u-\log p)^{q-j}\d u,\cr W_q(n,p)&:=\sum_{ 1\leqslant j\leqslant   q/2}{q\choose j}N_{j,q}(n,p).}\eqdef{defNjq}$$
An averaging argument over the largest prime factor of integers  enables one to essentially exploit the above in the form
 $$M_q(np)\leqslant 2M_q(n)+C\sum_{1\leqslant j\leqslant q/2}{q\choose j}M_j(n)M_{q-j}(n)\quad(p\nmid n)\eqdef{recMq2}$$
for some absolute constant $C$. This is the basis of various induction processes, leading to corresponding upper bounds. However, it turns out that an early appeal to Hšlder's inequality in the form
$$M_j(n)\leqslant M_q(n)^{(j-1)/(q-1)}\tau(n)^{(q-j)/(q-1)}\quad(n\geqslant 1,\,q\geqslant 2,\,1\leqslant j\leqslant q-1)$$
can only lead to upper bounds of the shape \eqref{GT85}, with or without the triple logarithm. The decisive progress achieved in \citer{KT23} consists in finding a framework inside which the full strength of an inequality like \eqref{recMq2} can be exploited. 
\par 
To prove \eqref{estsD}, we follow the approach of \citer{KT23} and insert a refinement  based on an appeal to Fourier transforms: see \eqref{Pars}. As noticed in \S\thinspace\ref{prgE*T}, this also yields an improvement on an estimate established in~\citer{HT88}---see \eqref{Tn0-div} {\it infra}.
\par
Since the theory of the $\Delta$-function already provided a number of surprises, it is rather risky to make a guess regarding the best possible exponent of $\log_2 x$ in \eqref{estsD}.  Regarding upper bounds, our approach might lead to the exponent 2, but we haven't been able to substantiate this so far.
  \par 
 It is worthwhile mentioning that the upper bound method displayed here (as of course that of \citer{KT23}) adapts smoothly to weighted averages of the type
 $$\sum_{n\leqslant x}g(n)\Delta(n)$$
 where  $g$ is a non-negative multiplicative function satisfying mild growth assumptions such as~\hbox{\citeplus{BT23}{(1.1),(1.2)}} and having mean-value 1 over the primes, viz.
 $$\sum_{p\leqslant x}g(p)\log p=x+O\big(x\e^{-(\log x)^c}\big)$$
 for some fixed $c>0$. This is precisely what is needed to handle averages of the $\Delta$-function at polynomial arguments. Thus we can state that if $F\in\z[X]$ is an irreducible polynomial with integer coefficients, then
 $$\sum_{n\leqslant x}\Delta\big(|F(n)|\big)\ll_F x(\log_2x)^{5/2}\qquad (x\geqslant 3).\eqdef{somDF}$$
  This is likely to be useful in applications, in particular to Waring's type problems: see \citer{BT25}.
  \par \medskip
  To conclude this introduction we note that the question of the quadratic mean-value of the $\Delta$-function has been solved, up to a power of $\log_2x$, in \citeplus{HT88}{Exercice 72}.
  We state the corresponding estimates in the following statement and, for the reader's convenience, provide the short proof in section \ref{prgDelta2}.
  \propt{D2}{ (\citer{HT88})}{We have
  $$x\log x\ll\sum_{n\leqslant x}\Delta(n)^2\ll x(\log x)(\log_2x)^2\qquad (x\geqslant  3).\eqdef{estD2}$$}
 A challenging open question is to determine the true exponent of $\log_2x$ in this problem.   
\goodbreak\medskip\bigskip
\paraun{Proof of \ref{thD}}
\paradeuxb{The set $E_T$}\drefdeux{2.1}
As in \citer{KT23}, the first step consists in defining a set of integers with useful multiplicative constraints.
\par 
Given a parameter $T\geqslant 3$,  put $$\eqalign{&\gb:={1\over \log 4-1},\quad f_T(y):=\delta\min\Big(\log (T\log 3y),{(\log_23y-\gb\log T)^2\over \log T}\Big)\quad(\delta>0,\,y\geqslant 1),\cr
&E:=\{n\geqslant 1:\mu(n)^2=1\},\quad n_y:=\prod_{p|n,\,p<y}p\quad(n\geqslant 1,\,y\geqslant 1),\cr
& P^+(n):=\max_{p|n}p\quad(n>1),\ P^+(1):=1,\cr
&E_{T}:=\big\{n\in E:\tau(n_y)\leqslant T(\log 3y)\e^{-f_T(y)} \ (y\geqslant 1)\big\}.\cr
& \cr}$$  
Here and throughout, $\mu$ denotes the Mšbius function, the letter $p$ denotes a prime number and~$\delta$ stands for an absolute, sufficiently small constant. It is useful   to bear in mind that $n\in E_T$ implies that $d\in E_T$ for all divisors $d$ of $n$.
\par \goodbreak
Given $x>y\geqslant 2$, let $\PP_{y,x}$ denote the probability on $E$ defined~by $$\PP_{y,x}(\{n\}):={1\over n}\prod_{y\leqslant  p<x}\Big(1+{1\over p}\Big)^{-1}\asymp{\log y\over n\log x}\qquad \big(P^-(n)\geqslant y,\,P^+(n)<x\big),$$
and let $\EE_{y,x}$ denote the corresponding expectation. Throughout we let $\som{y,x}$ indicate a summation over squarefree integers whose prime factors  are restricted to the interval $[y,x[$. When $y=2$, we simply write $\PP_x$, $\EE_x$, and $\som{x}$.
\par 
 It is proved in \citeplus{KT23}{prop. 5.1} that, provided $\delta$ is suitably chosen, we have
$$\PP_x( E\sset E_{T}) \ll {1\over T}\qquad (x\geqslant 2,\,T\geqslant 3).\eqdef{EssetET}$$

\par 
\bigskip\goodbreak
\paradeuxb{The set $E_T^{*}$}\drefdeux{E*T}
It is shown in \citer{KT23} that one can assume $M_2(n)/\tau(n)\leqslant T\log T$ $(n\in E_T)$ without perturbing \eqref{EssetET} and we note that this readily follows from \citeplus{HT88}{th. 46}, stated in terms of the function
$$T(n,0):=\sum_{\di{d,d'|n}{|\log (d'/d)|\leqslant 1}}1,$$
which satisfies $M_2(n)\leqslant T(n,0)\leqslant  4M_2(n)$ for all integers $n\geqslant 1$. From \citeplus{HT88}{thms. 46 \& 47} we actually have
$${1\over T\sqrt{\log T}}\ll\PP_x\big(M_2/\tau>T\big)\ll {\log T\over T}\qquad (T\geqslant 3).\eqdef{Tn0-div}$$
The next statement improves both bounds to
$1/T$. The upper estimate will be employed as a property of the~set
$$E_T^{*}:=\big\{n\in E_T:M_2(n)/\tau(n)\leqslant T\big\}.\eqdef{defET*}$$
It is useful to bear in mind that, since $M_2(n)/\tau(n)$ is multiplicatively increasing,\note{This follows from~\eqref{recMq}. See also \citeplus{HT88}{th. 40} for the variant related to $T(n,0)$.} we have $M_2(n_y)/\tau(n_y)\leqslant T$ for all $y\geqslant 1$ and all $n\in E_T^{*}$.

\Propp{q=2}{We have $$\PP_x( E\sset E_{T}^{*}) \ll {1\over T}\qquad (x\geqslant 2,\,T\geqslant 3).\eqdef{EssetET**}$$
Moreover
$$\PP_x\big(M_2/\tau>T\big)\asymp{1\over T}\cdot\eqdef{evalM2/tau}$$}
Since only \eqref{EssetET**} is necessary for the upper bound of \eqref{estsD}, we restrict below to the proof of this  estimate and postpone the proof of \eqref{evalM2/tau} to section \ref{prgM2T}.
\par \medskip
\noi{\it Proof of \eqref{EssetET**}.} By \eqref{EssetET}, it suffices to prove that $\PP_x(E_T\sset E_T^*)\ll1/T$. Put $\tau(n,\vartheta):=\sum_{d|n}d^{i\vartheta}$ $(n\geqslant 1,\vartheta\in\r)$. 
In view of \citeplus{HT88}{(3.2)}, Parseval's formula provides
$$M_2(n)={1\over 2\pi}\int_\r\Big({\sin\vartheta/2\over \vartheta/2}\Big)^2|\tau(n,\vartheta)|^2\d\vartheta \qquad (n\geqslant 1),\eqdef{ParsM2}$$
from which we derive, by applying the  Montgomery--Wirsing lemma  on mean values of Dirichlet polynomials (see, \eg, \citeplus{Te15}{lemma III.4.10}) to intervals $k-1/2\leqslant \vartheta\leqslant k+1/2$ $(k\in\z)$ and taking symmetry into account, that $${M_2(n)\over \tau(n)}\ll \int_0^{1}{|\tau(n,\vartheta)|^2\over \tau(n)}\d\vartheta\qquad (n\geqslant 1).\eqdef{Pars}$$
\goodbreak
The trivial bound $|\tau(n,\vartheta)|\leqslant \tau(n)$ immediately yields that the contribution of the range $0\leqslant \vartheta\leqslant 1/\log x$ has bounded expectation with respect to $\PP_x$.  Fix a large constant $K$, {and first assume that $T^K\leqslant \log x$.} Let  {then} $m_2(n,T)$ denote contribution  of the range $1/\log x\leqslant \vartheta\leqslant 1/T^K$ to the above integral, and let $M_2(n,T)$ denote the complementary contribution. Thus, in view of \eqref{EssetET}, we have, for a suitable absolute constant $c>0$, 
$$\PP_x\big(M_2/\tau>T\big)\ll {1\over T}+\alpha+\beta,$$
with 
$$\alpha:=\PP_x\Big(n\in E_T,\,m_2(n,T)>  c T\Big),\quad \beta:=\PP_x\Big(n\in E_T,\,M_2(n,T)>  cT\Big).$$
\par 
Under conditions $n\in  E_T$ and $K\geqslant K(\delta)$, we have $\tau(n_{\exp(1/\vartheta)})\leqslant  {(T/\vartheta)^{1-\delta}}$ whenever $\vartheta\leqslant 1/T^K$.  Therefore
$$\eqalign{\EE_x \big(m_2(n,T)\big)&\leqslant  \int_0^{1/T^K} \EE_x\bigg({ T^{1-\delta}|\tau(n,\vartheta)|^2\over\tau(n)\tau(n_{\exp(1/\vartheta)})\vartheta^{1-\delta}}\bigg) \d\vartheta\cr&\ll { T^{1-\delta}\over \log x}\int_0^{1/T^K}\som{x}_{n\geqslant 1}{|\tau(n,\vartheta)|^2\over n\tau(n)\tau(n_{\exp(1/\vartheta)})\vartheta^{1-\delta}}\dd\vartheta\ll  T^{1-\delta}\int_0^{1/T^K} {\dd\vartheta\over \vartheta^{1-\delta}},\cr}$$
where the inner sum, written as a product over primes, has been estimated by \citeplus{HT88}{lemma~30.1} or \citeplus{Te15}{lemma III.4.13}. For $K>1/\delta$, the above bound is $\ll1$, hence $\alpha\ll 1/T$.
\par 
Consider next
$$V_2(x):=\EE_x\big(M_2(n,T)^2\big)\ll\int_{1/T^K}^1\int_{1/T^K}^1F_x(\vartheta,\varphi)\d\vartheta\d\varphi,\eqdef{majV2}$$
with
$$\eqalign{F_x(\vartheta,\varphi)&:={1\over \log x}\som{x}_{n\in E_T}{|\tau(n,\vartheta)|^2|\tau(n,\varphi)|^2\over n\tau(n)^2}
\cdot\cr}$$
We use distinct upper bounds for this quantity, according to the sizes of $\tau(n_{\exp(1/\vartheta})$ and $\tau(n_{\exp(1/\varphi)})$. 
 By symmetry, we may assume $0< \vartheta\leqslant \varphi\leqslant 1.$ 
 Recall that for $n\in E_T$ we have 
 $$\tau(n_{\exp(1/\psi)})\leqslant {T\e^{-f_T(\exp(1/\psi))}\over \psi }\qquad \big(\psi\in\{\varphi,\vartheta\}\big).$$
\par For $j\geqslant 0$, we define
$$D_{j}(\varphi;T):=\bigg\{n\in E_T:{1\over 2^{j+1}}<{\varphi\tau(n_{\exp(1/\varphi)})\over T\e^{-f_T(\exp(1/\varphi))}}\leqslant  {1\over 2^j}
\bigg\},$$ 
let $F^{(j)}_x(\vartheta,\varphi)$ denote the subsum of $F_x(\vartheta,\varphi)$ restricted to summands $n$ from $D_{j}(\varphi;T)$ and let $V_2^{(j)}(x)$ denote the integral of $F^{(j)}_x(\vartheta,\varphi)$ over $[1/T^K,1]^2$. 
 Writing $\psi_p:=\psi\log p$ $(\psi>0)$, we have
 $$\leqalignno{U_1&:=\sum_{p\leqslant \exp(1/\varphi)}{1\over p}\leqslant \log\Big({1\over \varphi}\Big)+O(1),&\eqdef{estU1}\cr
U_2&:=\sum_{\exp(1/\varphi)<p\leqslant \exp(1/\vartheta)}{1+\cos\varphi_p\over p}\leqslant \log \Big({\varphi\over \vartheta}\Big)+O(1),&\eqdef{estU2}\cr
U_3&:=\sum_{\exp(1/\vartheta)<p\leqslant x}{(1+\cos\vartheta_p)(1+\cos\varphi_p)\over p}
\cr &\leqslant \log (1+\vartheta\log x)+\dm\log \Big({\varphi\log x\over 1+(\varphi-\vartheta)\log x}\Big)+O(1),
&\raise 4mm\hbox{\eqdef{estU3}}\cr}$$
where we used the formulae
$$\eqalign{(1+\cos\vartheta_p)&(1+\cos\varphi_p)\cr&=1+\cos\vartheta_p+\cos\varphi_p+\dm\cos(\varphi_p-\vartheta_p)+\dm\cos(\varphi_p+\vartheta_p),\cr
\sum_{p\leqslant y}{\cos \psi_p\over p}&=\log \Big({\log y\over 1+\psi\log y}\Big)+O(1)\quad(y\geqslant 2,\,0\leqslant \psi\leqslant 1).}\eqdef{sump}$$
\par \goodbreak
Let us now estimate $V_2^{(j)}(x)$. We have
$$F_x^{(j)}(\vartheta,\varphi)\ll {TG_x^{(j)}(\vartheta,\varphi)\e^{-f_T(\exp(1/\vartheta))/2-f_T(\exp(1/\varphi))/2}\over \sqrt{\vartheta\varphi}2^{j/2}\log x},$$
with $$G_x^{(j)}(\vartheta,\varphi):=\som{x}_{n\in D_{j}(\vartheta,\varphi;T)} {|\tau(n,\vartheta)|^2|\tau(n,\varphi)|^2\over n\tau(n)^2\sqrt{\tau(n_{\exp(1/\vartheta)})\tau(n_{\exp(1/\varphi)})}}\cdot$$
Observe that $n\in D_j(\vartheta,\varphi;T)$ implies
 $$ \omega(n_{\exp(1/\varphi)})=r_j(\varphi):=\fl{{\log ( {T/2^j\varphi})-f_T(\exp(1/\varphi))\over \log 2}}.$$ Split $n=uvw$, where $u:=n_{\exp(1/\varphi)}$ and where the prime factors of $v$, $w$, belong respectively to $\big]\exp(1/\varphi),\exp(1/\vartheta)\big]$, and $\big]\exp(1/\vartheta),x]$. We have
$${|\tau(n,\vartheta)|^2|\tau(n,\varphi)|^2\over n\tau(n)^2\sqrt{\tau(n_{\exp(1/\vartheta)})\tau(n_{\exp(1/\varphi)})}}\leqslant{\tau(u)|\tau(v, \varphi)|^2|\tau(w,\vartheta)|^2|\tau(w,\varphi)|^2\over uvw\sqrt{\tau(v)}\tau(w)^2}\cdot$$
Now, by \eqref{estU1},
$$\eqalign{\sum_{\di{u,v,w}{\omega(u)=r_j(\varphi)}} {\tau(u)|\tau(v, \varphi)|^2|\tau(w,\vartheta)|^2|\tau(w,\varphi)|^2\over uvw\sqrt{\tau(v)}\tau(w)^2}&\ll  { {(2U_1)}^{r_j(\varphi)}\e^{ {\sqrt{2}}U_2+U_3}\over r_j(\varphi)!}
\cr&\ll{\e^{2U_1+{\sqrt{2}}U_2 +U_3}\over \sqrt{1+U_1}}\ll{\e^{ {\sqrt{2}}U_2 +U_3} \over \varphi^2\sqrt{\log (3/\varphi)}} \cdot\cr}
 $$ Applying \eqref{estU2} and \eqref{estU3}, we infer that 
 $$\eqalign{{G_x^{(j)}(\vartheta,\varphi)\over \sqrt{\vartheta\varphi}\log x}&\ll{({\varphi/\vartheta})^{\sqrt{2}}(1+\vartheta\log x)\over \vartheta^{1/2}\varphi^{5/2}\sqrt{\log (3/\varphi) }\log x}\sqrt{\varphi\log x\over 1+(\varphi-\vartheta)\log x}\cr&\ll{\vartheta+1/\log x\over \varphi^{2-\sqrt{2}}\vartheta^{\sqrt{2}+1/2}\sqrt{\varphi-\vartheta+1/\log x}  \sqrt{\log ( 3/\varphi) }}\cdot\cr}$$
Ignoring the factor $\e^{-f_T}$ and splitting at $\varphi=2\vartheta$, we hence obtain that the  integral of $\varphi\mapsto F_x^{(j)}(\vartheta,\varphi)$  over $\vartheta\leqslant \varphi\leqslant 1$ is
$$\ll {T(\vartheta+1/\log x)^{3/2}\e^{-f_T(\exp(1/\vartheta))/2}\over 2^{j/2}\vartheta^{5/2}\sqrt{ \log ( 3/\vartheta)}} \cdot$$
Taking into account that $\vartheta\geqslant 1/\log x$ and integrating over the range $1/T^K\leqslant \vartheta\leqslant 1$ yields, exploiting here the factor $\e^{-f_T(\exp(1/\vartheta))}$, that $$V_2^{(j)}(x)\ll T/2^{j/2}.$$
Summing over $j\geqslant 0$, we obtain
$$V_2 (x)\ll  T.$$
\par 
Hence $\beta \ll V_2(x)/T^2\ll 1/T$, as required, provided $T^K\leqslant \log x$. However, if $T^K>\log x$, we have $\PP_x(M_2/\tau>T)\ll 1/T+\beta$ with now $\beta:=\PP_x\big(n\in E_T,\,M_2(n,(\log x)^{1/K})>  cT\big)$, and so the desired conclusion persists in view of the above calculations. \qed\goodbreak

\bigskip
\paradeuxb{Bounding moments inductively}\drefdeux{moms}
This is a simple reappraisal of the argument of \citer{KT23} which we reproduce with some simplifications  and some improvement arising from the  definition of the set $E_T^{*}$.\par 
Let  $\{\vartheta_{j,T}\}_{j=0}^{\infty}\in(\r^+)^{\N}$ be a sequence such that $ \vartheta_{0,T}= \vartheta_{1,T}=1$ and define the subset $E_{q,T}$ comprising those integers $n\in E_{T}^{*}$ such that 
$$M_j(n)\leqslant \tau(n)\vartheta_{j,T}\qquad (1\leqslant j\leqslant q),\leqno{(E_{q,T})}$$
with $ {E_{0,T}}=E_{1,T}=E_{T}^{*}$. Note that, by \eqref{recMq}, if $n\in E_{q,T}$, then any divisor of $n$ also lies in~$E_{q,T}$. By  
\ref{q=2}, provided $\vartheta_{2,T}\geqslant T$, we still have $$\PP_x\big(E\sset E_{2,T}\big)\ll 1/T\qquad (x\geqslant 2,\,T\geqslant 3).$$
\par \goodbreak
Recall notation $\som x$ and put 
$$S_q(x):=\PP_x\big(E_{q-1,T}\big)\EE_x\Big({M_q\over \tau}\Big|E_{q-1,T}\Big)=\prod_{p<x}\Big({1\over 1+1/p}\Big)\som{x}_{n\in E_{q-1,T}}{M_q(n)\over n\tau(n)}\cdot$$

\Propp{Sq}{Let $C_0$ be a large absolute constant. Assume $\vartheta_{1,T}\ge 1$, $\vartheta_{2,T}\ge C_0T$,  $\vartheta_{3,T}\geqslant 2C_0^2 T^2\log T$, and 
$$\vartheta_{j,T}\geqslant j! (C_0T)^{j-1}\quad(j\geqslant 2),\qquad \sum_{ 1\leqslant j\leqslant   q/2}{q\choose j}\vartheta_{j,T}\vartheta_{q-j,T}\leqslant {\vartheta_{q,T} \over C_0T\sqrt{\log T}}\quad(q\geqslant 3)\cdot\eqdef{inegthqT}$$
Then $$S_q(x)\leqslant {C_0\vartheta_{q,T}\over q^2T}\qquad (q\geqslant  3,\,x\geqslant 2,\,T\geqslant 3).\eqdef{majSq}$$}
\goodbreak\nid
Put $L_x:=\prod_{p<x}(1+1/p)$. Writing $n=mp$ with $P(m)<p$, we get
$$L_xS_q(x)\leqslant 1+\sum_{p<x} \som p_{m\in E_{q-1,T}}{M_q(pm)\over pm\tau(pm)}\cdot$$
Using $\tau(pm)= 2\tau(m)$ and expanding $M_q(pm)$ by \eqref{recMq} yields
$$L_xS_q(x)\leqslant 1+\sum_{p<x}\som p_{m\in E_{q-1,T}} {M_q(m)+W_q(m,p)\over pm\tau(m)}=\sum_{p<x} {L_pS_q(p)\over p }+V_q(x),$$
with
$$V_q(x):=1+\sum_{p<x} \som p_{m\in E_{q-1,T}}{W_q(m,p)\over pm\tau(m)}\cdot$$
Iterating this as in \citer{KT23}, we get
$$L_xS_q(x)\leqslant V_q(x)+\sum_{\di{2\leqslant P^+(n)<x}{n\in E }}{V_q(P^-(n))\over n}\ll 
V_q(x)+\sum_{p<x} {V_q(p)\over p}{\log x\over \log p},$$
hence
$$S_q(x)\ll {V_q(x)\over \log x}+\sum_{p<x}{V_q(p)\over p\log p}\cdot\eqdef{eqfonctT}$$
Indeed, at the first step we obtain
$$\sum_{p<x}{L_pS_q(p)\over p}\leqslant \sum_{p_1<x}{1\over p_1}\bigg\{\sum_{p_2<p_1}{L_{p_2}S_q(p_2)\over p_2}+V_q(p_1)\bigg\},$$
 then
$$\eqalign{\sum_{p<x}{L_pS_q(p)\over p}&\leqslant \sum_{p_1<x}{1\over p_1}\bigg(V_q(p_1)+\sum_{p_2<p_1}{1\over p_2}\Big\{\sum_{p_3<p_2}{L_{p_3}S_q(p_3)\over p_3}+V_q(p_2)\Big\}\bigg)\cr
&\leqslant \sum_{p_1<x}{V_q(p_1)\over p_1}+\sum_{p_2<p_1<x}{V_q(p_2)\over p_1p_2}+\sum_{p_3<p_2<p_1<x}{L_{p_3}S_q(p_3)\over p_1p_2p_3},\cr}$$
and so on.
\par Since $n_y=n_p$ when  {$P^+(n)<p<y\leqslant p^2$}, we may write,  still following \citer{KT23},
$$\eqalign{V_q(x)&\ll 1+\sum_{p<x}\som p_{m\in E_{q-1,T}}{W_q(m,p)\over pm\tau(m)} \int_p^{p^2}{\d y\over y\log 3y}\cr
&=1+\sum_{ 1\leqslant j\leqslant   q/2}{q\choose j}\int_1^{x^2}\som y_{m\in E_{q-1,T}}{1\over m\tau(m)}\sum_{\sqrt{y}\leqslant p< y}{N_{j,q}(m,p)\over p}{\d y\over y \log 3y },
\cr}$$ where $N_{j,q}$ is as defined in \eqref{defNjq}.
 By the argument leading to \citeplus{KT23}{(6.9)}, we get
$$V_q(x)\ll 1+\sum_{ 1\leqslant j\leqslant   q/2}{q\choose j}
\int_1^{x^2}\som y_{m\in E_{q-1,T}}\Big(1+2^j{\log y\over y^{1/4}}\Big){M_{q-j}(m)M_j(m)\over m\tau(m)}{\d y\over y(\log 3y)^2}\cdot$$
\goodbreak
Inverting summation and integration in \eqref{eqfonctT} and appealing to the uniform bound
$$\sum_{\sqrt{y}\leqslant p<x}{1\over p\log p}\ll  {1\over \log y},$$ we obtain
$$S_q(x)\ll 1+\sum_{ 1\leqslant j\leqslant   q/2}{q\choose j}
\int_1^{x^2}\som y_{m\in E_{q-1,T}}\Big(1+2^j{\log y\over y^{1/4}}\Big){M_{q-j}(m)M_j(m)\over m\tau(m)}{\d y\over y {(\log 3y)^3} }\cdot$$
\par 
The contribution of the term involving $y^{-1/4}$ may be handled as in \citer{KT23} using 
the trivial bound\note{Note that $M_j(m)\leqslant \Delta(m)^{j-1}M_1(m)\leqslant \tau(m)^j$ for $j\geqslant 1,\,m\geqslant 1$.}
$$M_{q-j}(m)M_j(m)\leqslant \tau(m)^{q}\leqslant (T\log y)^{q-2}\tau(m)^2,$$ valid whenever $m\in E_T$ and $P^+(m)< y$. This provides an overall term at most
$ C^{q}q!T^{q-2}$ where $C$ is a suitable absolute constant.
By the first assumption \eqref{inegthqT} we may choose $C_0$ such that   $$C^{q}q!T^{q-2}\leqslant C^2(C/C_0)^{q-2} q \vartheta_{q-1,T}\ll \vartheta_{q-1,T}/2^q .$$
\par 
Therefore
$$S_q(x)\ll {\vartheta_{q-1,T}\over 2^q}+\sum_{ 1\leqslant j\leqslant   q/2}{q\choose j}
\int_1^{x^2}{Z_{q,j}(y,T)\over y {(\log 3y)^3} }\d y,\eqdef{rec}$$
with 
$$Z_{q,j}(y,T):=\som y_{m\in E_{q-1,T}}{M_{q-j}(m)M_j(m)\over m\tau(m)}\cdot
$$
\par \goodbreak
We aim at establishing \eqref{majSq} by induction on $q\geqslant 3$. 
 By $(E_{q,T})$ and \eqref{defET*} we have, whenever $m$ is counted in  $Z_{q,j}(y,T)$,
$$M_j(m)\ll \vartheta_{j,T} T(\log 3y)\e^{-f_T(y)}\qquad (m\in E_{q-1,T},\,P^+(m)<y,\,1\leqslant j\leqslant q-1,\,y\geqslant 1).$$ 
Therefore
$$Z_{q,j}(y,T)\ll \vartheta_{j,T} T(\log 3y)^2S_{q-j}(y) {\e^{-f_T(y)}}.\eqdef{rec1}$$
\par  
When $q=3$, we infer from  $(E_{q,T})$ with $q=2$ and  the definition of  $E_T^*$ that
$$\eqalign{Z_{3,1}(y,T)&=\som y_{m\in E_{2,T}}{M_2(m)\over m }\leqslant\som y_{m\geqslant 1}{\min\big(T\tau(m),T^{3/2}\sqrt{\tau(m)\log y}\,\big)\over m }
\cr
&\ll{T^{3/2}(\log y)^2\over \sqrt{T}+(\log y)^{(3-\sqrt{8})/2}} \cdot\cr}\eqdef{majZ31}$$
Hence \eqref{rec} implies $$S_3(  x)\ll T\log T\ll   {\vartheta_{3,T}\over T}\cdot$$ This initiates the induction provided $C_0$ is sufficiently large.

If $q\geqslant 4$, then $q-j>2$ (except for $q=4$ and $j=2$), so we can use the induction hypothesis  \eqref{majSq}.  For $q=4$ and $j=2$,   we simply note that $Z_{4,2}(y,T)\leqslant  Z_{4,1}(y,T)$ 
since $M_2(m)^2\leqslant M_3(m)\tau(m)$ for all $m\geqslant 1$.  
 Appealing to \eqref{majSq} for bounding $S_{q-j}(y)$ for $1\leqslant j\leqslant q/2$,  
 we derive from \eqref{rec} and \eqref{rec1}
$$\eqalign{S_q(x)&\ll\sum_{1\leqslant j\leqslant q/2}{q\choose j}{\vartheta_{q-j,T}\vartheta_{j,T}\over q^2T}\int_1^{x^2}{T\e^{-f_T(y)}\over y(\log 3y)}\d y +{\vartheta_{q-1,T}\over 2^q}
\cr&\ll \sum_{1\leqslant j\leqslant q/2}{q\choose j}{\vartheta_{q-j,T}\vartheta_{j,T}\over q^2}\sqrt{\log T}.
}\eqdef{majSqavecsomme}$$
\par 
Under  assumption \eqref{inegthqT}, this implies \eqref{majSq}.
\qed\goodbreak
\bigskip
  
\paradeuxb{Completion of the argument}
 Observe that the sequence defined by 
$$\vartheta_{q,T}={q!\over q^2}\Big({2\pi^2 C_0\over 3}\Big)^{q-1} T^{q-1}(\log T)^{(q-1)/2}\quad(q\geqslant 1)$$  satisfies \eqref{inegthqT}.
\par 
Now, we put $\lambda:=C_1T (\log T)^{1/2}$  and aim at establishing
$$\PP_x\big(\Delta(n)>\lambda\log_2x\big)\ll{1\over T}\cdot\eqdef{majPxD}$$
By 
\eqref{EssetET**}, it is sufficient to show that 
$$\PP_x\big(n\in E_T^{*},\,\Delta(n)>\lambda\log_2x\big)\ll{1\over T}\cdot$$
\par 
However, from \eqref{majSq},  for $q\geqslant 3$,
$$\PP_x\big(E_{q-1,T}\sset E_{q,T}\big)\leqslant \PP_x\big(E_{q-1,T}\big)\EE_x\Big({M_q(n)\over \vartheta_{q,T}\tau(n)}\Big|E_{q-1,T}\Big)={S_q(x)\over \vartheta_{q-1,T}}\ll{1\over q^2T},$$
whence
$$\PP_x(E_T\sset \cap_{q\geqslant 3}E_{q,T})\ll{1\over T}\cdot$$
\par 
\goodbreak
Now observe that if $n\in\cap_{q\geqslant 3}E_{q,T}$ and $P^+(n)\leqslant x$, we have
$$\Delta(n)\leqslant 2 M_q(n)^{1/q}\ll \tau(n)^{1/q}\vartheta_{q,T}^{1/q}\ll (T\log x)^{1/q}qT^{1-1/q}(\log T)^{1/2},$$
where the first inequality is \citeplus{HT88}{(5.56)}.
Selecting $q=\fl{\log_2x}$ we have $\Delta(n)\leqslant \lambda\log_2x$ for a suitable large constant $C_1$ and hence~\eqref{majPxD}.
 \par 
It remains to estimate $\EE_x(\Delta).$  We may assume trivially that $10C_1\log_2x\leqslant \Delta(n)\leqslant (\log x)^{3}$. Indeed, 
$$\EE_x (\Delta, \Delta>(\log x)^{3})\leqslant  {\EE_x (\tau^2)\over (\log x)^3}\ll 1.$$   The contribution of those integers $n$ such that $2^j<\Delta(n)/\log_2x\leqslant 2^{j+1}$ is plainly
$$\ll{ 2^j j^{1/2}\log_2x\over 
2^j}= j^{1/2}\log_2x .$$
Summing over $j\ll \log_2x$, 
we finally get
$$\EE_x(\Delta )\ll(\log_2x)^{5/2}.$$
Since, by \citeplus{HT88}{th. 61}, we have
$S(x)\ll x\EE_x(\Delta),$ we obtain the upper bound of \eqref{estsD} as required.
\medskip\medskip
\paraun{Proof of  estimate  \eqref{evalM2/tau}}\drefun{M2T}
For integer $k$, let $\chi_k$ denote the indicator of the set $\big\{n\in E:\omega(n)=  k\big\}$. For notational convenience, we shall also write $\xi:=\sqrt{\log_2x}$ in the sequel. Recall the definition of the constant $\gb$ at the beginning of \S\thinspace\ref{prg2.1}
\Propl{Pom}{Let $0<\varepsilon<\dm$ be fixed. For integer $h$, $(\log x)^\varepsilon\leqslant 2^h\leqslant (\log x)^{1-\varepsilon}$, $y:=\exp2^h$, $\varepsilon\log_2x\leqslant k\leqslant (\log_2x)/\varepsilon$, we have
$$\EE_y\big(\chi_k\big)\asymp{(h\log 2)^k\over 2^hk!}\cdot\eqdef{Phm}$$
In particular, for $t\in h/\gb+[-\xi,\xi]$, we have
$$\EE_y(\chi_{h+t})\asymp {1\over 2^t\xi}\cdot\eqdef{Pht}$$}
\nid This readily follows from, \eg,  \citeplus{HR66}{lemma III.13}, which provides the estimate  $$\sum_{p_{1}<\cdots<p_{k}\leqslant y}{1\over p_{1}\cdots p_{k}}\asymp{(h\log 2)^k\over k!},$$
and Stirling's formula.\qed
\medskip
We may now embark on proving \eqref{evalM2/tau}. In view of \eqref{EssetET} and \eqref{EssetET**}, only the lower bound is to be established. Given $T\geqslant 3$ and a  large integer~$m$, put $$\xi_T:=\sqrt{\log T}, \quad\HH_T:=\Big({\gb\xi_T^2\over \log 2}+[0,\xi_T]\Big)\cap m\N,\quad \Y_T:=\{y_h:=\exp2^h:\,h\in\HH_T\}.$$ 
Let $K\geqslant 1$ to be chosen later and put $\sE_K:=\{n\in E:\max_{y\in\Y_T}n_y/y^K\leqslant 1\}$. We have
$${M_2(n)\over \tau(n)}\geqslant {M_2(n_y)\over \tau(n_y)}> {\tau(n_y)\over 2K\log y}\quad(n\in \sE_K,\, y\in\Y_T).$$ 
Let us then put $$t:={\log 2T\over  \log 2},\quad\gq_T(n):=\max_{y\in\Y_T}\chi_{h+t}(n_y),$$
so that $${M_2(n)\over \tau(n)}>{T\gq_T(n)\over 2K}\quad(n\in \sE_K).$$
\par \goodbreak
Consider
$\gd_T:=\EE_x\big(\gq_T\big).$ 
In view of \ref{Pom}, we first have
$$\EE_x\big(\chi_{h+t}(n_y)\big)\asymp{1\over T\xi_T},$$
hence
$$\gd_T\leqslant  \sum_{h\in\HH_T}\EE_x\big(\chi_{h+t}\big)\ll{1\over T}\cdot$$
Now, the inclusion-exclusion principle implies the lower bound
$$\gd_T\geqslant \sum_{h\in\HH_T}\EE_x\big(\chi_{h+t}\big)-\sum_{\di{h,j\in\HH_T}{h<j}}\EE_x\big(\chi_{h+t}\chi_{j+t}\big).$$
When $j=h+\ell$, $\ell\geqslant m$, we have
$$\EE_x\big(\chi_{h+t}\chi_{j+t}\big)\asymp {(\ell\log 2)^\ell\over T\xi_T2^\ell\,\ell!}\ll{\e^{-\beta\ell}\over T\xi_T},$$
with $\beta:=\log (2/\e\log 2)\approx0.05966$.
It follows that, for suitable absolute constants $c_1,\,c_2>0$,
$$\gd_T\geqslant {c_1\over mT}-{c_2\over T}\sum_{\di{\ell\geqslant m}{\ell\md0m}}{\e^{-\beta\ell}}\geqslant {c_1/m-c_2/(\e^{\beta m}-1)\over T}\cdot$$
Selecting $m$ sufficiently large, we get that $\gd_T\gg 1/T$.
\par 
To conclude the proof it only remains to observe that
$$\eqalign{\EE_x\Big(\1_{E\sset \sE_K}\gq_T\Big)&\ll{1\over \e^K}\sum_{h\in\HH_T}\EE_x\Big(\chi_{h+t}(n_y)n_y^{1/\log y}\Big)\cr&\ll{1\over \e^K}\sum_{h\in\HH_T}{1\over 2^{h}(h+t)!}\Big(\sum_{p\leqslant y}{1\over p^{1-1/\log y}}\Big)^{h+t}\ll{1\over \e^K T},\cr}$$
and get, with a suitable choice of constant $K$,
$$\PP_x\Big(M_2(n)/\tau(n)>T/2K\Big)\gg 1/T.$$
 Up to substituting $2KT$ to $T$, this implies the lower bound in \eqref{evalM2/tau}.
\smallskip 
\medskip
 \paraun{Proof of {\ref{D2}}}\drefun{Delta2}
 The lower bound readily follows from the pigeon-hole principle in the form  $$\Delta(n)\gg\tau(n)/\log 2n\quad(n\geqslant 1).$$ \par 
 In order to establish the upper bound, we first note that by a standard argument (see, \eg, \citeplus{Te85}{(22)}, \citeplus{HT88}{th. 61}), it is sufficient to show that 
 $$\EE_x(\Delta^2)\ll(\log_2x)^2\log x.\eqdef{D2log}$$ Now, we
 observe that, for all integers $r\geqslant 1$, $s\geqslant 1$ and suitable $u\in \r$, we have
 $$\eqalign{\Delta(rs)&\leqslant \sum_{\di{d|r,\,t|s}{-1/2\leqslant \log dt-u\leqslant 1/2}}1\leqslant 2\sum_{d|r}\sum_{t|s}(1-|\log dt-u|)^+\cr
 &={1\over \pi}\int_\r\Big({\sin\vartheta/2\over \vartheta/2}\Big)^2\tau(r,\vartheta)\tau(s,\vartheta)\e^{-iu\vartheta}\d\vartheta.}$$
 By \eqref{ParsM2} and the Cauchy-Schwarz inequality, we get as in \citeplus{HT88}{Exercise 47}
$$\Delta(rs)^2\leqslant 4 M_2(r)M_2(s)\qquad (r³1,\,s\geqslant 1).$$
\goodbreak
Therefore
$$\mu(n)^2\Delta(n)^2\leqslant \sum_{rs=n}{\mu(r)^2M_2(r)\over \tau(r)}{\mu(s)^2M_2(s)\over \tau(s)}\qquad (n\geqslant 1).$$
It is then sufficient to apply \eqref{Pars} and \citeplus{HT88}{th. 42}, which readily implies that $M_2(n)\mu(n)^2/\tau(n)$ has average order $\ll\log_2n$, to get \eqref{D2log}, as required.  
\bigskip
\noi{\bf Acknowledgement.} The authors take pleasure in addressing warm thanks to Terence Tao, for pointing out an inaccuracy in a preliminary version of this article and to Dimitris Koukoulopoulos for pointing out several in  successive drafts. 
\goodbreak
\bigskip  
 \centerline{\twelvebf References}
\bigskip
 {\eightpoint\leftskip9mm\rightskip5mm
 \bibtem{BT23} R. de la Bretche \& G. Tenenbaum, Two upper bounds for the Erd\H os--Hooley Delta-function, {\it Sci. China Math. \bf66}, no. 12 (2023), 2683Ð2692.\par 
 \bibtem{BT25} R. de la Bretche \& G. Tenenbaum, Mean values of arithmetic functions
and application to sums of powers, https://arxiv.org/abs/2403.19320v2.\par 
 \bibtem{Er73} P. Erd\H os, Problem 218, Solution by the proposer, {\it Can. Math. Bull. \bf 16} (1973),  621-622.\par 
 \bibtem{EN75} P. Erd\H os \& J.-L. Nicolas, RŽpartition des nombres superabondants, {\it Bull. Soc. math. France \bf 103} (1975), 65-90.\par 
 \bibtem{FGK23} K. Ford,  B. Green, D. Koukoulopoulos, Equal sums in random sets and the concentration of divisors, 
{\it Inventiones Math.} {\bf 232}, (2023),1027--1160. \par 
\bibtem{FKT23} K. Ford, D. Koukoulopoulos, T. Tao, A lower bound on the mean value of the Erd\H os--Hooley Delta function, {\it Proc. London Math. Soc.}, to appear.
\par 
\bibtem{HR66} H. Halberstam \& K.F. Roth,
 {\it Sequences}, Oxford, 1966,  2nd. ed. Springer, 1983.
\par 
\bibtem{HT82} R.R. Hall \& G. Tenenbaum, On the average and normal orders of Hooley's $\Delta$-function,
{\it J. London Math. Soc.} (2) {\bf 25} (1982), 392-406. 
\bibtem{HT84} R.R. Hall \& G. Tenenbaum, The average orders of Hooley's $\Delta_r$-functions, {\it
Mathematika \bf 31} (1984), 98--109.
 \bibtem{HT86} R.R. Hall \& G. Tenenbaum, The average orders of Hooley's $\Delta_r$-functions, II,
 {\it Compositio Math. \bf 60} (1986), 163--186.\par 
  \bibtem{KT23} D. Koukoulopoulos \& T. Tao, An upper bound on the mean value of the Erd\H os--Hooley Delta function, {\it Proc. London Math. Soc. (3) \bf 127}, no. 6 (2023),  1865Ð1885.
 \bibtem{HT88} R.R. Hall \& G. Tenenbaum, {\it Divisors}, Cambridge tracts in
mathematics, no 90, Cambridge University Press (1988).\par 
\bibtem{Ho79} C. Hooley, On a new technique and its applications to the theory of numbers, {\it Proc.
London Math. Soc.} (3) {\bf 38} (1979), 115--151.
\par
\bibtem{Te85} G. Tenenbaum,
Sur la concentration moyenne des diviseurs, {\it Comment. Math. Helvetici} {\bf 60} (1985), 411-428. 
\bibtem{Te86} G. Tenenbaum,
Fonctions $\Delta$ de Hooley et applications, SŽminaire de thŽorie des nombres, Paris 1984-85, Prog. Math. 63 (1986), 225-239.
 \par
\bibtem{Te15} G. Tenenbaum, {\it Introduction to analytic and probabilistic number theory}, 3rd ed., Graduate Studies in Mathematics 163, Amer. Math. Soc. (2015).\par
\bibtem{Te16}  G. Tenenbaum, Fonctions multiplicatives, sommes d'exponentielles, et loi des grands nombres, {\it Indag. Math. \bf27} (2016), 590-600.
\par }
{\leftskip2mm\rightskip-2cm\sevenrm
\gutter=4cm \doublecolumns
 \obeylines \baselineskip=7pt
RŽgis de la Bretche
UniversitŽ Paris CitŽ, Sorbonne UniversitŽ, CNRS,
Institut Universitaire de France,
Institut de Math. de Jussieu-Paris Rive Gauche
 F-75013 Paris  
France
\smallskip
{\seventt regis.de-la-breteche@imj-prg.fr}
\phantom a\goodbreak
 GŽrald Tenenbaum
Institut \'Elie Cartan
Universit\'e de Lorraine
 BP 70239\par
54506 Vand\oe uvre-ls-Nancy Cedex
 France
\smallskip
{\seventt gerald.tenenbaum@univ-lorraine.fr}
\singlecolumn
\par}
\end